\documentclass{amsart}

\numberwithin{equation}{section}

\usepackage{amssymb}

\newtheorem{thmm}{Theorem}
\newtheorem{thm}{Theorem}[section]
\newtheorem{lem}[thm]{Lemma}

\newtheorem{rem}[thm]{Remark}

\newcommand\B{{\mathcal B}}
\newcommand\Co{{\mathcal C}}
\newcommand\D{{\mathcal D}}

\newcommand\Or{{\mathcal O}}

\newcommand\C{{\mathbb C}}

\newcommand\N{{\mathbb N}}
\newcommand\R{{\mathbb R}}

\newcommand\ve{\varepsilon}

\newcommand\vf{\varphi}

\newcommand\de{d^\flat}
\newcommand\supp{\operatorname{supp}}

\newcommand\Id{\text{\bf Id}}
\newcommand\Tr{\text{Tr\,}}
\newcommand\Exp{\operatorname{Exp}}

\begin{document}

\title{Zeta functions and Dynamical Systems}
\author{Carlangelo Liverani}
\author{Masato Tsujii}
\address{Carlangelo Liverani\\
Dipartimento di Matematica\\
II Universit\`{a} di Roma (Tor Vergata)\\
Via della Ricerca Scientifica, 00133 Roma, Italy.}
\email{{\tt liverani@mat.uniroma2.it}}
\address{Masato Tsujii\\ Department of Mathematics\\ Hokkaido University\\ Sapporo, 060-0810, Japan }
\email{{\tt tsujii@math.sci.hokudai.ac.jp}}
\date{\today}
\thanks{The first named author would like to thank D.Dolgopyat for
pointing out to him the wonders of tensor products.}
\begin{abstract}
In this brief note we present a very simple strategy to investigate
dynamical determinants for uniformly hyperbolic systems. The construction builds
on the recent introduction of suitable functional spaces which allow
to transform simple heuristic arguments in rigorous ones. Although the
results so obtained are not exactly optimal the straightforwardness of
the argument makes it noticeable.
\end{abstract}
\maketitle

\section{introduction}
The goal of the paper is to investigate the properties of the dynamical
Fredholm determinants of uniformly hyperbolic systems and to relate
them to the statistical properties of such systems. This subject has
been widely investigated and there exists a large literature where
many partial results are obtained. We refer the reader to \cite{Babook} for
references and an introduction to the subject, to \cite{Li,BT2} for a more
recent account of the situation and to \cite{chaos} for an in depth
discussion of the physical relevance of these issues.

The basic idea presented in this paper is to study the action of the dynamics on an
appropriate singular functional kernel, as suggested in Dmitry Dolgopyat's
thesis (Princeton 1997), to obtain results on the radius of convergence of the dynamical
Fredholm determinant, its relation to the spectral properties of the
transfer operator and the Ruelle resonances. The new ingredient
allowing to carry out such a program is the possibility, after \cite{GL} and
\cite{BT}, to introduce spaces in which such singular kernels are legal object.
To clarify matters we start with a folklore explanation.

Let $X$ be a $d$-dimensional $\Co^{r+1}$ Riemannian manifold and $T:X\to X$ a
$\Co^{r+1}$ diffeomorphism which satisfies some hyperbolicity
condition. (We assume at least that all the periodic points of $T$ are
hyperbolic.) For each $g\in\Co^{r}(X,\C)$ we
define the Ruelle transfer operator $T_g:\Co^{r}(X,\C)\to \Co^{r}(X,\C)$
by
\begin{equation}
\label{eq:dual-pf}
T_g h:= g\cdot h\circ T.
\end{equation}
The dynamical Fredholm determinant of this operator $T_g$ is formally
defined by\footnote{Note that if the number of periodic points does
not grow more than exponentially, then $d_{T,g}^{\flat}(z)$ is well
defined and
holomorphic in a sufficiently small disk.}
\begin{equation}
\label{eq:fredholm}
d_{T,g}^{\flat}(z)=\exp \left[- \sum_{n\ge 1}\frac{z^n}{n} \sum_{x\in\text{\rm Fix
}T^{n}}\frac{g_n(x)}{|\det (Id-DT^n(x))|}\right],
\end{equation}
where $g_n(x):=\prod_{i=0}^{n-1} g(T^i(x))$.

Let us first note that this can be heuristically regarded as the determinant $\det
(Id-z \cdot T_g)$. Indeed, let $\delta$ be the
distribution on $X^2=X\times X$ defined by
\[
\delta(h)=\int_{X} h(x,x) dx.
\]
Then the ``kernel'' of the operator $T_g^n$ is given by $
(Id\otimes T_g)(\delta)$.\footnote{Formally, $\delta$ is $\delta(x-y)$
(where now $\delta$ is the physicists {\sl delta function})
and the kernel is given by $g(x)\delta(T_gx-y)$. As we shall see later,
the action of $Id\otimes
T_g$ can in fact be extended to an operator on the space of distributions.}
Thus, as in the case of operators with smooth kernel, it would be
natural to define
\[
\Tr T_g =\langle\delta,(Id\otimes T_g)(\delta)\rangle.
\]
Though the product of two distribution is not defined in general, we
will be able to give an appropriate meaning to the right hand side above
since the singular supports (in $\Co^r$ sense) of $\delta$ and
$(Id\otimes T_g)(\delta)$ do not intersect. We find \footnote{One
can easily guess this formula by approximating $\delta$ by a sequence
of $\Co^{\infty}$ functions, see \cite{Babook} or \cite{Li} if details
are really needed. }
\[
\langle\delta,(Id\otimes T_g)^n(\delta)\rangle=
\sum_{x\in\text{\rm Fix
}T^{n}}\frac{g_n(x)}{|\det (Id-DT^n(x))|}.
\]
The definition \eqref{eq:fredholm} of the dynamical Fredholm
determinant follows then via the formal relation $\det(Id -zA)=\exp
(-\sum_{n=1}^{\infty}(z^n \Tr A^n)/n)$.

It is thus natural to expect that the properties of the dynamical
Fredholm determinant as holomorphic function are closely related to
the spectral properties of the operator $T_g$. In this paper, we
present an argument providing exactly such a relation, although in a
slightly more restrictive setting. The argument is rigorous, yet it
follows the above simple ideas very closely.

Let $T:X\to X$ be an Anosov diffeomorphism, i.e. there
exists a $DT$-invariant decomposition $TM=E^u\oplus E^s$ and constants
$\lambda\in(0,1)$ and $C>0$ such that:
\[
\|DT^n|_{E^s}\|\le C\lambda^n,\qquad
\|DT^{-n}|_{E^u}\|\le C\lambda^n \qquad\mbox{for all \ $n\ge 0$.}
\]
In \cite{GL} and \cite{BT}, Banach spaces $\B$ of distributions on $X$
are defined so that the operator $T_g$  extends to a bounded operator
$T_g:\B\to \B$ whose essential spectral radius is bounded by
$\|g\|_{L^{\infty}}\cdot \lambda^{\alpha_r}$, where
$\alpha_r:=\min\{[r/2],r-[r/2]\}$, $[a]\in\N$ being the closest
integer to $a\in\R$.\footnote{Actually, \cite{BT} allows the better
bound $\alpha_r=r/2$. 
Also, \cite{GL} deals
explicitly only with the adjoint of $T_g$ in the
case $g\equiv 1$ (SRB measures), yet the extension to the present setting is
straightforward.
}
In addition, it is shown that
the eigenvalues outside the essential spectral radius have a well
defined dynamical meaning (Ruelle resonances). Let
$\rho_*=\|g\|_{L^{\infty}}\cdot \lambda^{\alpha_r/2}$. Our result is
as follows:

\begin{thmm}\label{th}
 $d_{T,g}^{\flat}(z)$ extends holomorphically to
$D(\rho_*^{-1})=\{|z|<\rho_*^{-1}\}$ and the zeros of such an
extension are in one-one correspondence, with multiplicity, to the
inverse of the eigenvalues of $T_g:\B\to \B$ in the region $\{|z|>\rho_*\}$.
\end{thmm}

This result is not new nor optimal. Kitaev \cite{Ki} has given a stronger result
for the extendibility part of the former claim,\footnote{Essentially,
instead of the bound $\rho_*^{-1}$, Kitaev has the
more natural bound  $(\|g\|_{L^{\infty}}\cdot \lambda^{r/2})^{-1}\sim
(\|g\|_{L^{\infty}}\cdot \lambda^{\alpha_r})^{-1}$,
that is the inverse of the bound for the essential spectral radius of
$T_g:\B\to \B$.} while the spectral
interpretation, albeit for a smaller radius, appeared already in
\cite{Li} and, more recently, in \cite{BT2} it has been obtained for a domain corresponding to
the result of Kitaev. Nevertheless, the proofs yielding such sharper
results are far more complex than the present argument.

The non-optimality of the above theorem is the price for
considering, in the following, the operator $T_g^*\otimes T_g$ instead
of the operator $Id\otimes T_g$ used in the previous heuristic
argument. Unfortunately, we do not
know how to treat the transfer operator $ Id\otimes T_g$ directly as
the mapping $Id\times T$ is not hyperbolic and the extension of the
results in \cite{GL}, \cite{BT} to the partially hyperbolic setting
is far from trivial (if possible at all).

\section{Basic definitions}
For $g\in\Co^{r}(X,\C)$ we define the transfer operator
$T_g:\Co^{r}(X,\C)\to \Co^{r}(X,\C)$ by
\[
T_g h:= g \cdot (h\circ T).
\]
Its formal 
adjoint (in fact,
dual)
$T_g^*:\Co^{r}(X,\C)\circlearrowleft$ is given by
the transfer operator
\[
T_g^* f=|\det DT^{-1}|\cdot (g\cdot f)\circ T^{-1}.
\]
That is, $\langle\overline{ T_g h}, f\rangle_{L^2(X)}=\langle \bar h, T_g^* f\rangle_{L^2(X)}$.

Let $\D'_{r}(X)$ be the space of distribution on $X$ of order
$r$.\footnote{Here, by $\D'_r(X)$ we mean
the dual of the space $\Co^r(X)$ defined as follows. For $r\geq 0$,
let $\lfloor r\rfloor$ be its integer part. We denote
by $\bar \Co^r$ the set of functions which are $\lfloor r\rfloor$
times continuously differentiable, and whose $\lfloor r\rfloor$-th
derivative is H\"{o}lder continuous of exponent $r-\lfloor r\rfloor$
if $r$ is not an integer. To fix notation, in this paper we choose,
for each $r\in\R_+$, a norm on $\bar\Co^r$ functions so that $|\vf_1
\vf_2|_{\Co^r} \leq |\vf_1|_{\Co^r} |\vf_2|_{\Co^r}$. We will denote
by $\Co^r$ the closure in $\bar \Co^r$ of the set of $\Co^\infty$
functions. It coincides with $\bar \Co^r$ if $r$ is an integer, but
is strictly included in it otherwise. In any case, it contains $\bar
\Co^{r'}$ for all $r'>r$.}
Using the above formal relation, we can extends the  operators
$T_g$ and $T_g^*$ to continuous operators $T_g:\D'_{r}(X)\to
\D'_{r}(X)$ and $T_g^*:\D'_{r}(X)\to
\D'_{r}(X)$, respectively.

Next we define
$T_g^*\otimes T_g:\Co^r(X^2,\C)\to \Co^r(X^2,\C)$ as the unique
extension of  $T_g^*\otimes
T_g:\Co^r(X,\C)\otimes\Co^r(X,\C)\circlearrowleft$.
The latter operator reads
\begin{equation}\label{eq:tensor}
T_g^*\otimes T_g(\vf)(x,y):=g(T^{-1}x)\cdot |\det D_xT^{-1}|\cdot g(y)\cdot \vf(T^{-1}x,Ty).
\end{equation}
So it can be interpreted as the transfer operator associated to the hyperbolic
mapping $T^{-1}\times T$ with the
weight $\widetilde{g}(x,y)=g(T^{-1}x)\cdot |\det D_xT^{-1}|\cdot g(y)$. Note that the formal
adjoint of $T_g^*\otimes T_g$ is $T_g\otimes T^*_g$. As above, we can
extend these operators to
\begin{equation}
\label{eq:d-transfer}
T_g\otimes T_g^*(\vf): \D'_{r}(X^2)\to \D'_{r}(X^2)
\quad\mbox{and}\quad
T^*_g\otimes T_g(\vf): \D'_{r}(X^2)\to \D'_{r}(X^2).
\end{equation}

\section{The Banach spaces of distribution}
The spaces of distributions on which we have defined the transfer
operators are not appropriate for a study of the dynamics.
Many recent
works (e.g. \cite{GL,BT,FR}) have focussed on the problem of finding adapted functional
spaces and different choices have different advantages. Accordingly, it may
better not to focus on a particular choice but to enlightened which are
the properties needed to carry out the study of the Zeta functions. We
will list the properties of the function spaces that suffice for our
argument. Yet, for definiteness of exposition we will comment explicitly the
Banach spaces introduced in \cite{GL} and remark that they indeed satisfy such properties.

In \cite{GL}, S.Gou\"ezel and the first-named-author introduced a
scale of Banach spaces $\B^{p,q}$ with $q\in \R_+$, $p\in \N$ and
$p+q<r$ adapted to $\Co^{r+1}$ Anosov diffeomorphisms $T:X\to X$. The parameters
$p$ and $q$ will be fixed at the end of the argument. We
denote $\B=\B^{p,q}$ and set
\[
\rho=\rho_{p,q}=\lambda^{\min\{p,q\}} \|g\|_{L^{\infty}},\quad \widetilde\rho=\rho \|g\|_{L^{\infty}}.
\]
The basic properties of $\B$ are the following (see \cite{GL} for a proof):
\begin{itemize}
\item[{\bf (P1)}]  $\Co^{r}(X)$ is continuously embedded in $\B$ and
its image is a dense subset.
\item[{\bf (P2)}]
$\B$ is continuously embedded in $\D'_r(X)$.
\item[{\bf (P3)}]
$T_g:\B\to\B$ is a bounded operator with essential spectral radius bounded by $\rho$.
\end{itemize}
\begin{rem}
The meaning of {\bf (P3)} is that we use {\bf (P2)} to identify $\B$
with a subspace of $\D'_r(X)$ and consider the restriction of $T_g$ to
$\B$. With such an identification the embedding in {\bf (P1)} is
required to be the standard embedding of $\Co^\infty$ in $\D'_r(X)$.
In the following we will use the embeddings {\bf (P2)}, (and {\bf (P5)}) to identify
the elements of $\B$  (and $\widetilde \B$) with distributions without
making further remarks.
\end{rem}
As already noted in \eqref{eq:tensor} the operator $T_g^*\otimes T_g$ is a transfer operator
for the Anosov diffeomorphism $T^{-1}\times T:X^2\to X^2$ with the
same hyperbolicity constant $\lambda$ of $T$. So, as above, we can
introduce a Banach space $\widetilde{\B}$ with the following
properties
\begin{itemize}
\item[{\bf (P4)}]  $\Co^{r}(X^2)$ is continuously embedded in
$\widetilde{\B}$ and its image is a dense subset.
\item[{\bf (P5)}]
$\widetilde{\B}$ is continuously embedded in $\D'_r(X^2)$.
\item[{\bf (P6)}]
$T_g^*\otimes T_g:\widetilde{\B}\to\widetilde{\B}$ is a bounded
operator with essential spectral radius bounded by $\widetilde\rho$.
\end{itemize}
The reader should be aware that usually there is some freedom in the definition of the Banach
spaces. For example, in \cite{GL} they depend on a family
$\Sigma$ of admissible leafs, that is, $\Co^{r+1}$ embedded compact
$\dim E^s$ dimensional submanifolds with boundary close to local
stable manifolds. By taking the family appropriately,\footnote{In the
definition of $\widetilde{\B}$, one can take $\Sigma$ so that the diagonal in
$X^2$ is covered by finitely many elements in $\Sigma$. This
immediately implies {\bf (P8)}.} we can insure that $\widetilde{\B}$
enjoys the following extra properties
\begin{itemize}
\item[{\bf (P7)}] $(T_g^*\otimes T_g)^{n_0}(\delta)$ is  contained in
$\widetilde{\B}$ for some $n_0\in\N$.

\item[{\bf (P8)}] The functional $\bar\delta:\Co^{\infty}(X^2)\to \C$,
$\bar\delta(\varphi)=\delta(\varphi)$, extends to a bounded functional
$\bar\delta:\widetilde{\B}\to \C$.
\end{itemize}

Finally, as it should be apparent from the previous heuristic
argument, we need some control on how to approximate singular kernels
by smooth ones.
Let $\{(U_i,\Psi_i:U\to \R^n)\}_{i=1}^k$ be a $\Co^{r+1}$ atlas of $X$, and let
$\{\rho_i\}_{i=1}^k$ be a $\Co^\infty$ partition of unity subordinated to such an atlas. Next,
define the functions $j_\ve\in\Co^\infty(\R^d,\R_+)$ so that
$\int_{\R^d}j_\ve=1$ and $\supp(j_\ve)\subset\{x\in\R^d\;:\;\|x\|\leq
\ve\}$. We then define \footnote{Note that $J_\ve$ is
well defined for $\ve$ small enough.}
\[
J_\ve f(x):=\sum_i\int_{\R^d}\rho_i\circ\Psi_i^{-1}(y)j_\ve(\Psi_i(x)-y)
f\circ \Psi_i^{-1}(y) dy\quad\mbox{for $f\in \Co^r(X)$.}
\]
Let $J^*_\ve$ be the formal adjoint of $J_{\ve}$. Clearly
these extend to  bounded operators
\[
J_{\ve}:\D'_r(X)\to\Co ^{\infty}(X)\quad\mbox{and}\quad
J^*_{\ve}:\D'_r(X)\to \Co^{\infty}(X)\, .
\]
We also define
\[
\widetilde{J}_{\ve}:\D'_r(X^2)\to \Co^{\infty}(X^2)
\]
as the unique extension of $J_\ve\otimes J_\ve:\Co^{r}(X)\otimes \Co^r(X) \to \Co^{\infty}(X^2)$.
Then we have\footnote{These properties are essentially proven in \cite[section 7]{GL}, or see \cite{Li}.}
\begin{itemize}
\item[{\bf (P9)}] For $u\in \B$, $J_\ve^* J_\ve u\to u$ in $\B$,  as $\ve\to +0$.
\item[{\bf (P10)}] For $u\in \widetilde{\B}$, we have $\widetilde{J}_\ve u\to u$ in $\widetilde{\B}$,  as $\ve\to +0$.
\end{itemize}

\begin{rem}From now on we will use only the above properties, regardless of the
way the spaces are actually constructed.
\end{rem}
Here are two consequence of properties  ({\bf P1})-({\bf P10})
that show the relevance of the above objects for the problem at hand.
\begin{lem}\label{lem:deltabar}
Let $\{\ell_i\}\subset\B'$ and $\{e_i\}\subset\B$.
If $\sum_{i=1}^{k} \ell_i\otimes e_i$ belongs to
$\widetilde{\B}$, then we have
\[
\bar\delta\left(\sum_{i=1}^{k} \ell_i\otimes e_i\right)
=\sum_{i=1}^{k}\ell_i(e_i).
\]
\end{lem}
\begin{proof}
Define $J'_\ve:\B'\to\Co^\infty$ by the duality relation
$J'_\ve\ell(h):=\ell(J_\ve^* h)$.
By ({\bf P10}), $\widetilde{J}_{\ve}(\sum_{i=1}^{k} \ell_i\otimes e_i)=
\sum_{i=1}^{k} J'_\ve(\ell_i) \otimes J_{\ve}(e_i)\in \Co^{\infty}(X^2)$
converges to $\sum_{i=1}^{k} \ell_i\otimes e_i$ in $\widetilde{\B}$,
as $\ve\to +0$. Since
\[
\begin{split}
\bar\delta\left(\sum_{i=1}^{k} J'_\ve(\ell_i) \otimes
J_{\ve}(e_i)\right)&= \sum_{i=1}^{k} \langle \overline{J'_\ve(\ell_i)},
J_{\ve}(e_i)\rangle_{L^2(X)} =\sum_{i=1}^{k} J'_\ve \ell_i( J_{\ve}(e_i))\\
&=\sum_{i=1}^{k}  \ell_i(J_\ve^* J_{\ve}(e_i))
\end{split}
\]
the claim of the lemma holds by ({\bf P8}) and ({\bf P9}).
\end{proof}

\begin{lem}\label{lem:formulae}
Set $\delta_T=(T^*_g\otimes \Id) \delta$. For each $n\ge n_0$, we have
\[
\begin{split}
&\bar\delta((T_g^*\otimes T_g)^n\delta)=\sum_{x\in\text{\rm
Fix}\; T^{2n}}\frac{g_{2n}(x)}{|\det(\Id-DT^{2n}(x))|}\\
&\bar\delta((T_g^*\otimes T_g)^n\delta_T)=\sum_{x\in\text{\rm Fix
}T^{2n+1}}
\frac{g_{2n+1}(x)}{|\det(\Id-DT^{2n+1}(x))|}
\end{split}
\]
\end{lem}
\begin{proof}
Let $\delta_{\ve}=\widetilde{J}_{\ve}((T_g^*\otimes
T_g)^{n_0}\delta)$. By ({\bf P6}), ({\bf P7}) and ({\bf P10}) we have
$(T_g^*\otimes T_g)^{n-n_0}\delta_{\ve}\to (T_g^*\otimes T_g)^n\delta$ in $\widetilde{\B}$.
We leave it to the reader to check that  the right hand side of the
first equality above is obtained as the limit $\lim_{\ve\to +0}
\bar\delta((T_g^*\otimes T_g)^{n-n_0}\delta_{\ve})$.\footnote{If in trouble,
see \cite{Li} for details} The second equality is
obtained in a parallel manner.
\end{proof}

\section{The proof}

Take $\sigma>\max\{\rho,\widetilde\rho\}$ arbitrarily.  Then
properties ({\bf P3}) and ({\bf P6}) imply
\begin{equation}\label{decomp}
T_g^*=P+R:\B\to \B,\qquad
T_g^* \otimes T_g=\widetilde{P}+\widetilde{R}:\widetilde{\B}\to \widetilde{\B}
\end{equation}
where $P$ and $\widetilde P$ are of finite rank, the spectral radius
of $R$ and $\widetilde R$ are bounded by $\sigma$, and  $P R=RP=0$,
$\widetilde{P}\widetilde{R}=\widetilde{R}\widetilde{P}=0$.
Notice that, for each $h,f\in\Co^\infty(X)$, $n\in\N$,
\begin{equation}
\label{eq:powers}
\begin{split}
((T_g^*\otimes T_g)^n\delta)(h\otimes f)&
=\delta((T_g)^n h\otimes (T_g^*)^{n}f)= \langle
\overline{(T_g)^n h}, (T_g^*)^{n}f\rangle_{L^2(X)}\\
&=\langle \bar h, (T_g^*)^{2n}f\rangle_{L^2(X)}= (T_g^{2n}h)(f),
\end{split}
\end{equation}
where, in the last expression, $h$ is interpreted as an element of
$\D'_r$ (and hence $\B$) via the natural embedding.
It follows, using the Neumann series for $|z|$ small enough,
\begin{equation}
\label{eq:resolvent-rel}
(\Id-z(T_g^*\otimes T_g))^{-1}(T_g^*\otimes T_g)^{n_0}\delta(h\otimes f)
=((\Id-z(T_g)^2)^{-1}T_g^{2n_0}h)(f)
\end{equation}
where $n_0$ is the integer in the condition ({\bf P7}).
Note that, from (\ref{decomp}) and the arbitrariness of $\sigma$, both
sides of \eqref{eq:resolvent-rel}
have a meromorphic extensions to an closed disk
$\{z\in\C\;;\;|z|\leq \sigma^{-1}\}$,
and the equality must holds for those extensions.
Since $R^{2n}$ and $\widetilde{R}^n$ can be written as\footnote{Here one
has to choose $\sigma$ so that no eigenvalues belong to the circle
$\{z \in\C\;:\;|z|=\sigma\}$.}
\[
\begin{split}
R^{2n}&=\frac 1{2\pi i}\int_{|\zeta|=\sigma}z^n(z\cdot \Id
-(T_g)^2)^{-1} dz\;;
\\  \widetilde{R}^n&=\frac 1{2\pi
i}\int_{|\zeta|=\sigma}z^n(z\cdot \Id -(T_g^*\otimes T_g))^{-1}
dz\end{split}
\]
we have, for all $n> n_0$ and each $h,f\in\Co^\infty(X)$,
\begin{equation}\label{eq:finiterank}
\begin{split}
\widetilde{R}^n \delta(h\otimes
f)&:=\widetilde{R}^{n-n_0}(T_g^*\otimes T_g)^{n_0}\delta(h\otimes
f)=(R^{2n}(h))(f),\\
\widetilde{P}^{n} \delta(h\otimes f)&:=\widetilde{P}^{n-n_0}
(T_g^*\otimes T_g)^{n_0} \delta(h\otimes f)=(P^{2n}h)(f)
\end{split}
\end{equation}

Next, write the finite rank operator $P^{2n}$ as
$P^{2n}g=\sum_{i=1}^{k}\ell_i^{(n)}(g) e_i^{(n)} $ where
$\ell_i^{(n)}\in\B'$ and $e_i^{(n)}\in\B$. By \eqref{eq:finiterank},
for each $h,f\in\Co^\infty(X)$ and $n> n_0$,
\[
\widetilde{P}^n \delta(h\otimes
f)=\sum_{i=1}^{k}\ell_i^{(n)}(h)\cdot e_i^{(n)}(f)=\sum_{i=1}^{k}
\ell_i^{(n)}\otimes e_i^{(n)}(h\otimes f).
\]
Since
$\Co^\infty(X)\otimes\Co^\infty(X)$ is dense in $\Co^p(X^2)$ in the
$\Co^{p}$ topology,\footnote{This is a direct consequence of
Stone-Weierstrass theorem.} and since the elements of $\widetilde{\B}$ are
distributions of order $p$ by hypothesis,
it follows, for each $n> n_0$,
\begin{equation}
\label{eq:structure}
\widetilde\B\ni
\widetilde{P}^n \delta=\sum_{i=1}^{k}\ell_i^{(n)}\otimes e_i^{(n)}\in
\B'\times \B.
\end{equation}
We can thus apply Lemma \ref{lem:deltabar} obtaining
\[
\bar\delta\left(\widetilde{P}^n \delta\right) =
\sum_{i=1}^{k}\ell_i^{(n)}(e_i^{(n)})=\Tr P^{2n}.
\]
In conclusion we have, for $n> n_0$,
\begin{align*}
\bar\delta\left((T_g^*\otimes T_g)^n\delta\right)&=
\bar\delta\left( \widetilde{P}^n\delta\right)+
\bar\delta\left( \widetilde{R}^n\delta\right)
=\Tr P^{2n}+\Or(\sigma^{n}).
\end{align*}
By replacing $\delta$ by $\delta_T=(T_g^*\otimes \Id)\delta$ in the argument above, we obtain also
\begin{align*}
\bar\delta\left( (T_g^*\otimes T_g)^n\delta_T\right)&=
\bar\delta\left( \widetilde{P}^n\delta_T\right)
+
\bar\delta\left(\widetilde{R}^n\delta_T\right)
=\Tr P^{2n+1}+\Or(\sigma^{n}).\end{align*}
Using Lemma \ref{lem:formulae} we can conclude
\begin{equation}\label{eq:zetaf}
\begin{split}
\de_{T,g}(z):&=\Exp\left(-\sum_{n=1}^\infty\frac{z^{n}}n
\sum_{x\in\text{Fix }T^n}
\frac{g_n(x)}{|\det(\Id-D_{x}T^n)|}\right)\\
&=\Exp\left(-\sum_{n=1}^\infty\frac{z^{n}}n \Tr
P^{n}+p_0(z)+\sum_{n=1}^\infty\Or(\sigma^{n/2})z^n\right) \\
&=\det(\Id -zP)\Exp\left(-p_0(z)-\sum_{n=1}^\infty\Or(\sigma^{n/2})z^n\right) ,
\end{split}
\end{equation}
where $p_0(z)$ is a polynomial of order $2n_0+1$.

Since the last series is convergent for $|z|<\sigma^{-1/2}$, we have
that the $\de_{T,g}(z)$ is holomorphic in such a disk and the zeroes
correspond to the eigenvalues of $P$, that is to the eigenvalues of
$T_g:\B\to \B$ in the region $\{|z|>\sigma^{1/2}\}$.

Applying the argument above to the case where $p=[r/2]$ and $q$ is
arbitrarily close to $r-[r/2]$ we obtain Theorem 1. Let us conclude by
reiterating the generality of the approach.

\begin{rem}
Given any Banach spaces $\B$ and $\widetilde{\B}$ satisfying the properties {\bf
(P1-10)} the proof above applies, hence we obtain Theorem 1 with
$\rho_*=\max\{\rho^{-1},\widetilde{\rho}^{-1/2}\}$.
\end{rem}

\end{document}